\documentclass[reqno,12pt]{amsart} 
\usepackage[top=1.5in,right=1.125in,left=1.125in,bottom=1.5in]{geometry}
\usepackage{amssymb}
\usepackage{amsmath}
\usepackage{graphicx}
\usepackage{color}
\definecolor{red}{rgb}{0.7,0,0}
\definecolor{grey}{RGB}{112,112,112}
\definecolor{blue}{RGB}{034,113,179}
\usepackage[colorlinks=true,citecolor=blue,linkcolor=red,urlcolor=grey,backref=page]{hyperref}

\makeatletter
\@namedef{subjclassname@2020}{\textup{2020} Mathematics Subject Classification}
\makeatother

\newcommand{\koniec}{\begin{flushright}  $\Box $ \end{flushright}}
\newtheorem{theo}{Theorem}[section] 
\newtheorem{prop}[theo]{Proposition}  
\newtheorem{lemma}[theo]{Lemma}

\newtheorem{col}[theo]{Corollary}

\theoremstyle{remark}

\newtheorem{remark}[theo]{Remark}

\newcounter{mnotecount}[section]

\renewcommand{\themnotecount}{\thesection.\arabic{mnotecount}}

\newcommand{\mnote}[1]
{\protect{\stepcounter{mnotecount}}$^{\mbox{\footnotesize
$
\bullet$\themnotecount}}$ \marginpar{
\raggedright\tiny\em
$\!\!\!\!\!\!\,\bullet$\themnotecount: #1} }

\newcommand{\hook}{{\setlength{\unitlength}{11pt}   
                   \begin{picture}(.833,.8)
                   \put(.15,.08){\line(1,0){.35}}
                   \put(.5,.08){\line(0,1){.5}}
                   \end{picture}}}

\newcommand{\C}{\mathbb{C}}
\newcommand{\Z}{\mathbb{Z}}

\newcommand{\R}{\mathbb{R}}

\newcommand{\Rho}{\mathrm{P}}

\newcommand{\om}{\omega}

\def\p{\partial}
\def\be{\begin{equation}}

\def\ee{\end{equation}}

\def\bea{\begin{eqnarray}}
\def\eea{\end{eqnarray}}

\newcommand{\cT}{\mathcal{T}}

\numberwithin{equation}{section}
\begin{document} 


 \def\sideremark#1{\ifvmode\leavevmode\fi\vadjust{\vbox to0pt{\vss
  \hbox to 0pt{\hskip\hsize\hskip1em
  \vbox{\hsize2.6cm\tiny\raggedright\pretolerance10000
   \noindent #1\hfill}\hss}\vbox to8pt{\vfil}\vss}}}%
                        
                                                    %
 \newcommand{\edz}[1]{\sideremark{#1}}

\title{Conformally  K\"ahler structures}
\author{Maciej Dunajski}
\address{Department of Applied Mathematics and Theoretical Physics\\ 
  University of Cambridge\\ Wilberforce Road, Cambridge CB3 0WA, UK.\\
  and\\
  Faculty of Physics\\
  University of Warsaw
Pasteura 5, 02-093 Warsaw, Poland
}
\email{m.dunajski@damtp.cam.ac.uk}
\author{Rod Gover}
\address{Department of Mathematics\\
 The University of Auckland, \\
Private Bag 92019, Auckland 1142, New Zealand}
\email{r.gover@auckland.ac.nz}
\begin{abstract}
We establish a one-to-one correspondence between K\"ahler metrics in a
given conformal class and parallel sections of a certain vector bundle
with conformally invariant connection, where the parallel sections
satisfy a set of non--linear algebraic constraints that we
describe. The vector bundle captures 2-form prolongations and is isomorphic
to $\Lambda^3(\cT)$, where ${\cT}$ is the tractor bundle of conformal
geometry, but the resulting connection differs from the normal
tractor connection by curvature terms.

 Our analysis leads to a set of obstructions for a Riemannian metric
to be conformal to a K\"ahler metric. In particular we find an
 explicit algebraic condition for a Weyl tensor which must hold if
 there exists a conformal Killing-Yano tensor, which is a necessary
 condition for a metric to be conformal to K\"ahler. This gives an
 invariant characterisation of algebraically special Riemannian
 metrics of type $D$ in dimensions higher than four.
\end{abstract}


\maketitle
\begin{center}
{\em In memory of Anna Edmunds (1969–2023)}
\end{center}
\section{Introduction}
Let $(M, g)$ be a Riemannian manifold of even dimension $n\geq 4$. Does there exist a non--zero function $\Omega:M\longrightarrow
\R$ such that $\hat{g}=\Omega^2 g$ is K\"ahler?  That is does there exists a
non--degenerate two--form $\hat{\omega}$ which is covariantly constant
with respect to the Levi--Civita connection of $\hat{g}$ and such that
the endomorphism $J:TM\rightarrow TM$ defined\footnote{This is one of
  many equivalent definitions. It implies that $J$ is hermitian and
  that
\[
[T^{(1, 0)}, T^{(1, 0)}]\subset T^{(1, 0)},\quad 
\mbox{where}\quad T^{(1, 0)}=\{X\in TM\otimes\C, J(X)=iX\}.
\]}  by 
${\hat{\omega}(X, Y)}=\hat{g}(X, JY)$ satisfies $J^2=-\mbox{Id}$?

 In this paper we shall focus on local obstructions, which arise
 because the conformal to K\"ahler problem leads to an
 over--determined system of PDEs of finite type \cite{branson}. We
 shall establish a one--to--one correspondence between K\"ahler
 metrics in a conformal class and certain ({\it special}) parallel sections of a
 vector bundle $E\rightarrow M$ of rank $n(n+1)(n+2)/6$ \be
\label{bundle}
E\cong \left(\begin{array}{c}
\Lambda^2(M)\\ \Lambda^1(M)\oplus\Lambda^3(M)\\ \Lambda^2(M)
\end{array} \right)
\ee
equipped with a connection that takes the form
\be
\label{connection}
{\quad{\mathcal{D}} \left(\begin{array}{c}
\omega\\ 
K,\;\mu\\ 
\Sigma 
\end{array} \right)= 
\left(\begin{array}{c} \nabla\omega-{\mu}- g\oslash {K} \\ \nabla K  -{\Sigma}-P\oslash{\omega} {,}\;\;\;             \nabla\mu - g\oslash{\Sigma}
- P\oslash{\omega} -C\oslash {\omega}
\\ \nabla\Sigma-P\oslash {K}-A\oslash {\omega}- C\oslash {K}
\end{array} \right).}
\ee Here $A=$Cotton, $C=$Weyl and $P=$Schouten are different
components of the curvature tensor of $g$ and its derivatives, and
$\oslash$ indicates an algebraic operation involving contractions of
various kinds (differing in each line above). The latter will be
specified in Theorem \ref{theo01} in \S\ref{section2}.  We shall say
that a section $\Psi=(\omega, K, \mu, \Sigma)$ of $(E, {\mathcal D})$
is {\it special} if it satisfies a set of algebraic conditions \be
\label{constrains_s}
\mathcal{Q} (\Psi)=0,
\ee 
which will be specified in Theorem \ref{prop2} in \S\ref{section4}.
It is the presence of these conditions which makes the analysis difficult. 
General parallel sections of $(E, {\mathcal D})$ are in one--to--one correspondence
with conformal Killing--Yano tensors \cite{tachi,mason,GS}, and the algebraic constraints  single out 
those conformal Killing--Yano tensors which give rise to K\"ahler forms. 

The analysis leads to two kinds of obstructions. Those arising from reducing the holonomy
of the curvature of ${\mathcal D}$ to a subgroup stabilising a section of $E$,
and those arising from differentiating the algebraic conditions. This second class
of constraints is an overdetermined system of algebraic equations which can admit
non--zero solutions if and only if the relevant Bezout resultants 
vanish (Theorem \ref{theobezou}). The dimension of the variety of constraints
(\ref{constrains_s}) in the fibres of $E$ is at most $(n^2+2n+4)/4$ if $n>4$. In dimension $n=4$ the constraints can be solved
explicitly, reducing the rank $20$ vector bundle $E$ to a rank $10$ 
sub-bundle
\[
E_0=\Lambda^2_+(M)\oplus \Lambda^1(M)\oplus \Lambda^2_-(M),
\]
whose section consists of a self--dual two--form, a one--form, and an anti--self--dual 
two--form \cite{DT} (this is isomorphic to the bundle $\Lambda^3_+({\bf{T}})$ of self--dual
tractor three-forms). That is, using the notation above, the system of $\mathcal{Q} (\Psi)=0$ is equivalent to
\[
*\omega=\omega, \quad K=*\mu, \quad *\Sigma=-\Sigma.
\]
The complete set of obstructions has been constructed explicitly in this 
case\footnote{The details of the prolongation connection on $E_0$ constructed in \cite{DT} differ from
those in Theorem \ref{theo01} as in four dimensions it is convenient to fix the freedom in
the definition of $\Sigma$ so that $\Sigma$ is anti--self--dual and so $g(\omega, \Sigma)=0$.} \cite{DT}.
In \S\ref{sectractors} we shall link Theorem \ref{theo01} and Theorem \ref{prop2} with the tractor approach to conformal geometry
\cite{BEG}, and identify the prolongation bundle $E$ with the third exterior power of the rank $(n+2)$ tractor
bundle $\cT$ over $M$.  The prolongation connection (\ref{connection}) differs from the standard tractor connection by the curvature terms, and some of the non--linear constraints can be encoded in algebraic conditions involving the scale tractor.

\subsection*{Acknowledgements} Both authors acknowledge support
from the Royal Society of New Zealand via Marsden Grants 19-UOA-008 and 24-UOA-005.
MD is also grateful to the University of Auckland, and similarly RG to
the University of Cambridge, for the hospitality during visits when
this work was carried over. The authors would also like to thank the
Isaac Newton Institute for Mathematical Sciences, Cambridge, for
support and hospitality during the programme Twistor theory, where
work on this paper was completed. This work was supported by EPSRC
grant EP/Z000580/1. RG was also suppported by a Simons Foundation
Fellowship during this period.

\section{Prolongation of the conformal-to-K\"ahler system}
\label{section2}
In this section we shall directly construct the prolongation connection
(\ref{connection}) underling the conformal-to-K\"ahler problem.
\begin{theo}
\label{theo01}
There exists a correspondence between K\"ahler metrics in a given
conformal class and parallel sections $\Psi=(\omega, K, \mu, \Sigma)$
of $(E\rightarrow M, \mathcal D)$ given by (\ref{bundle}) and
(\ref{connection}).
\end{theo}
\noindent
{\bf Remark.} This is  not a one--to--one correspondence. Every K{\"a}hler metric corresponds to a parallel section of $(E, \mathcal{D})$,
but not all parallel sections give rise to K{\"a}hler metrics. Those which do will be characterised in Theorem \ref{prop2}.\\
{\bf Proof of Theorem {\ref{theo01}}.} Consider a K\"ahler structure $(\hat{g}, \hat{\omega})$, and define $(g, \omega)$ 
by\footnote{The standard scaling in Kahler geometry would be $\hat{\omega}=\Omega^2\omega$.
We attach a different conformal weight to $\omega$ to ensure that
equation (\ref{mdprol0}) is conformally invariant.}
\[
\hat{g}=\Omega^2g, \qquad
\hat{\omega}=\Omega^3 \omega ,
\]
where $\Omega$ is a smooth positive function.
The condition $\hat{\nabla}\hat{\omega}=0$ yields
\be
\label{mdprol0}
\nabla_{a}\omega_{bc}=\mu_{abc}+2g_{a[b}K_{c]},
\ee
where $\mu\in\Lambda^3(M)$ and $K\in\Lambda^1(M)$ are given by
\be
\label{solution_1}
\mu=-3\Upsilon\wedge\omega,\quad K=-\Upsilon\hook\omega,
\ee
and $\Upsilon =\Omega^{-1}d\Omega$. 

Conversely, assume that (\ref{mdprol0}) holds with arbitrary $\mu$ and $K$, for some $\omega$ such that ${\omega^a}_b{\omega^b}_c$ is pure trace. 
Recall the decomposition of the Riemann tensor in conformal geometry
\[
R_{abcd}=C_{abcd}+P_{ac}g_{bd}-P_{bc}g_{ad}+P_{bd}g_{ac}-P_{ad}g_{bc},
\]
where $C_{abcd}$ is the Weyl curvature and
\[
P_{ab}=\frac{1}{n-2}\Big(R_{ab}+\frac{R}{2(1-n)}g_{ab}\Big)
\]
is the  Schouten  tensor. Under conformal rescalings $\hat{g}=\Omega^2g$ we have
\[
\hat{C}_{abcd}=\Omega^2C_{abcd}, \quad \hat{P}_{ab}=P_{ab}-\nabla_a\Upsilon_b
+\Upsilon_a\Upsilon_b-\frac{1}{2}|\Upsilon|^2g_{ab}.
\]
We differentiate (\ref{mdprol0}),
commute the derivatives, and use the Ricci identity
\[
[\nabla_a, \nabla_b]\omega_{cd}=-{R_{abc}}^p\omega_{dp}
+{R_{abd}}^p\omega_{cp}.
\]
This leads to a set of 
algebraic conditions
\be
\label{weyl}
{C_{bc[a}}^e\omega_{d]e}+{C_{ad[b}}^e\omega_{c]e}=0,
\ee
and a pair of linear differential equations
\be
\label{mdprol1}
\nabla_a K_b={P_a}^c\omega_{bc}+\Sigma_{ab}
\ee
and
\be
\label{mdprol2}
\nabla_a \mu_{bcd}=-3g_{a[b}\Sigma_{cd]}-3 P_{a[b}\omega_{cd]}
-\frac{3}{2}{C_{[bc|a}}^p\omega_{p|d]},
\ee
where $\Sigma_{ab}$ is as an yet undetermined two--form. Differentiating 
(\ref{mdprol1}) and (\ref{mdprol2}) once
more gives
\be
\label{mdprol3}
\nabla_a\Sigma_{bc}=2P_{a[b}K_{c]}-{P_a}^e\mu_{ebc}+\frac{1}{2}{A^p}_{bc}\omega_{pa}
+{A^p}_{a[b}\omega_{c]p}+{C_{bca}}^pK_p,
\ee
where $A_{abc}=\nabla_bP_{ca}-\nabla_cP_{ba}$ is the Cotton tensor.
The system is now closed, as derivatives of all unknowns have been determined.
We can combine equations (\ref{mdprol0}), (\ref{mdprol1}),   
(\ref{mdprol2}), (\ref{mdprol3}) into a connection (\ref{connection}), where the 
meaning of $\oslash$ in each slot is now clear. \koniec
As a spin off from the prolongation procedure we deduce the 
following (well known)
\begin{col}
If a non K\"ahler manifold $(M, g)$ is Einstein and conformal to K\"ahler, then 
$g$ admits a Killing vector. 
\end{col}
\noindent
{\bf Proof.}
This follows directly from (\ref{mdprol1}).
If $g$ is Einstein then the RHS of (\ref{mdprol1}) is skew-symmetric, and thus $K$ satisfies the Killing equations $\nabla_{(a}K_{b)}=0$.\koniec

In four dimensions one can establish a stronger result:
an ASD Einstein metric with non--zero Ricci scalar is conformal to K\"ahler
if and only if it admits a Killing vector \cite{Derdzinski, DT}.
\section{Type $D$ and obstructions algebraic in Weyl tensor}
\label{section3}
If we view both  $C$ and $\omega$ as endomorphisms of $\Lambda^2$ given by
\[
C(\phi)_{ab}={C^{cd}}_{ab}\phi_{cd}\quad \mbox{and}\quad \omega(\phi)_{ab}={\omega_{[a}}^c\phi_{b]c},
\]
then the constraint (\ref{weyl}) is equivalent to the commutativity of these endomorphisms. 
Thus, they can be diagonalised in the same basis. In \cite{mason} it was used 
to show that the Weyl tensor is of algebraic type $D$ in the sense of \cite{coley, pravda}. 
\vskip5pt
We adopt a different approach. Consider a linear map 
\begin{eqnarray*}
&&B:\Lambda^2\rightarrow \mathcal{W}\subset\Lambda^2\odot\Lambda^2, \quad\mbox{with}\\
&& 
B(\phi)_{bcad}:= {C_{bc[a}}^e\phi_{d]e}+{C_{ad[b}}^e\phi_{c]e},
\end{eqnarray*}
given by the
LHS of equation (\ref{weyl}). Here $\mathcal{W}$ is a vector space of
rank-four tensors which have the algebraic symmetries of a Weyl
tensor, 
i.e. if $e\in \Gamma(\mathcal{W})$ then
\[
e_{a[bcd]}=0, \quad e_{abcd}=e_{[ab]cd}, \quad  e_{abcd}=e_{ab[cd]},
\]
and $e$ is  trace free with respect to metric contractions.
The dimension of $\mathcal{W}$ is greater than that of $\Lambda^2$, and
equation  (\ref{weyl})
implies that $B$ has a non--empty kernel which contains a non--degenerate two--form.
Therefore the rank of 
$B$ is not maximal. This leads to a set of algebraic conditions  on the Weyl tensor 
which we shall now give.
\begin{theo}
\label{theobezou}
Let $X\in \Lambda^{2}(TM)$ be a bi--vector and let $\beta_X:\Lambda^2\rightarrow
\Lambda^2$ be given by
\be
\label{defibeta}
{({\beta_X})_{ab}}^{cd}\equiv X^{ef}{\beta_{efab}}^{cd}, \quad\mbox{where}\quad
{\beta_{bcad}}^{ef}=
{C_{bc[a}}^e{\delta^f}_{d]}+{C_{ad[b}}^e{\delta^f}_{c]}.
\ee
Then $g$ is conformal to a K\"ahler metric only if, for all bi--vectors $X$, 
$\mbox{det}(\mathcal B)=0$,
where
\be
\label{b_matrix}
{\mathcal B}=\frac{1}{N!}
\left(\begin{array}{cccccc}
0 & (N-1)! & 0 & \dots  & 0 & 0 \\
s_2 & 0 & (N-2)! & \dots  & 0 &0\\
s_3 & s_2 & 0& (N-3)! & \dots   &0\\
\vdots &  \vdots &  & \ddots& \vdots & \vdots \\
s_{N-1} &  s_{N-2} &  & & 0 & 1 \\
s_N & s_{N-1} & s_{N-2} & \dots &s_2 & 0
\end{array}\right),
\ee
and $s_k\equiv \mbox{Tr}({\beta_X}^k)$, and $k=1, 2, \dots, N$.
\end{theo}
\noindent
{\bf Proof.}
Rewriting (\ref{weyl}) as
\[
{\beta_{bcad}}^{ef}\omega_{ef}=0,
\]
where $\beta$ is given by (\ref{defibeta}),
we deduce that for any fixed values of the pair of indices $[bc]$ the
determinant of the resulting $n(n-1)/2$ by $n(n-1)/2$ traceless square matrix $\beta$ must vanish. 
In the case $n=4$ this leads to an invariant condition on the self--dual part of Weyl
tensor, as explained in \cite{DT}.
For any bi--vector $X\equiv X^{ab}$ consider a composition of homomorphisms
\[
\textstyle\Lambda^2(M)\xrightarrow{\; B\;}{\mathcal W}
\xrightarrow{\;X\hook\underbar{\enskip}\;}
\Lambda^2(M)
\]
where the second map is a contraction. This gives a traceless
homomorphism
\[
{({\beta_X})_{ab}}^{cd}\equiv X^{ef}{\beta_{efab}}^{cd}.
\]
To find an invariant obstruction - a tensor of rank $n(n-1)$ on $TM$ - 
we shall use
the Cayley--Hamilton theorem for traceless $N\times N$ matrices, where
$N=n(n-1)/2$ is the dimension of $\Lambda^{2}(M)$. Set
\[
s_k=\mbox{Tr}({\beta_X}^k), \quad k=1, 2, \dots, N
\]
so that
\[
s_1=0, \quad s_2=X^{ab}X^{cd}  
{\beta_{abpq}}^{rs} {\beta_{cdrs}}^{pq}, \quad
s_3=X^{ab}X^{cd} X^{ef} 
{\beta_{abpq}}^{rs} {\beta_{cdrs}}^{uv} {\beta_{efuv}}^{pq}, \quad \dots.
\]
The determinant of $\beta_X$ can then be expressed as the $N$th Bell 
polynomial \cite{bell}
\begin{eqnarray}
\label{obstruction1}
\mbox{det}(\beta_X)&=&\frac{1}{N!}B_N(s_1, -1! s_2, 2! s_3, \dots, (-1)^N N! s_N)\\
&=&\mbox{det}(\mathcal B),\nonumber
\end{eqnarray}
where ${\mathcal B}$ is given by (\ref{b_matrix}).
\koniec
For this to be a non--trivial obstruction we need to show that $\beta_X$ can have maximal rank
(and therefore is injective) for some Weyl tensor ${C_{abc}}^d$ as it will then have maximal rank
in a neighbourhood of this Weyl tensor in ${\mathcal W}$. We could show it by specifying an element
of ${\mathcal W}$ at a point in $M$, but we can do even better, and write down a metric
which gives rise to such an injective Weyl tensor. On an open set in $\R^6$ with
coordinates $(x, y, z, t, u, v)$
consider a metric
\[
g=dx^2+dy^2+dz^2+dt^2+du^2+dv^2+c(t^2+yt)dxdy+ c(t^2+tu)dudv,
\]
where $c$ is a constant,
and take
\[
X=2\p_x\wedge\p_y+\p_z\wedge\p_u.
\]
Evaluating the Weyl tensor of this metric at the point
$(0, 0, 0, 0, 0, 0)$, and computing the obstruction
(\ref{obstruction1})   yields
\begin{eqnarray*}
\mathcal B
&=&
-{\frac {{\it s_3}\,{{\it s_4}}^{3}}{1152}}+{\frac {{{\it s_2}}^{6}{\it 
s_3}}{138240}}-{\frac {{{\it s_2}}^{3}{{\it s_3}}^{3}}{7776}}+{\frac {{
\it s_3}\,{{\it s_6}}^{2}}{216}}-{\frac {{{\it s_2}}^{5}{\it s_5}}{19200}}
-{\frac {{{\it s_3}}^{3}{\it s_6}}{972}}+{\frac {{{\it s_2}}^{4}{\it s_7}
}{2688}}+{\frac {{{\it s_4}}^{2}{\it s_7}}{224}}-{\frac {{\it s_7}\,{\it 
s_8}}{56}}\\
&&-{\frac {{{\it s_2}}^{3}{\it s_9}}{432}}+{\frac {{{\it s_3}}^{2}
{\it s_9}}{162}}-{\frac {{\it s_6}\,{\it s_9}}{54}}-{\frac {1}{50}}\,{
\it s_{10}}\,{\it s_5}+{\frac {{\it s_{11}}\,{{\it s_2}}^{2}}{88}}-\frac{1}{44}\,{\it 
s_{11}}\,{\it s_4}-\frac{1}{36}\,{\it s_3}\,{\it s_{12}}-{\frac {{{\it s_2}}^{2}{\it s_3
}\,{\it s_8}}{192}}+\\
&&
{\frac {{\it s_2}\,{\it s_5}\,{\it s_8}}{80}}+{\frac{
{\it s_3}\,{\it s_4}\,{\it s_8}}{96}}
-{\frac {{{\it s_2}}^{2}{\it s_4}\,{
\it s_7}}{224}}-{\frac {{\it s_2}\,{{\it s_3}}^{2}{\it s_7}}{252}}+
{\frac 
{{\it s_2}\,{\it s_6}\,{\it s_7}}{84}}+
{\frac {{\it s_3}\,{\it s_5}\,{\it 
s_7}}{105}}+\\
&&{\frac {{{\it s_2}}^{3}{\it s_3}\,{\it s_6}}{864}}-{\frac {{{
\it s_2}}^{2}{\it s_5}\,{\it s_6}}{240}}+{\frac {{\it s_4}\,{\it s_5}\,{
\it s_6}}{120}}+{\frac {{{\it s_2}}^{3}{\it s_4}\,{\it s_5}}{960}}+
{\frac 
{{{\it s_2}}^{2}{{\it s_3}}^{2}{\it s_5}}{720}}-{\frac {{\it s_2}\,{\it s_3
}\,{{\it s_5}}^{2}}{300}}\\
&&-{\frac {{\it s_2}\,{{\it s_4}}^{2}{\it s_5}}{320
}}
-{\frac {{{\it s_3}}^{2}{\it s_4}\,{\it s_5}}{360}}-{\frac {{{\it s_2}}^
{4}{\it s_3}\,{\it s_4}}{4608}}+
{\frac {{{\it s_2}}^{2}{\it s_3}\,{{\it s_4
}}^{2}}{768}}+{\frac {{\it s_2}\,{{\it s_3}}^{3}{\it s_4}}{1296}}-\\
&&
{\frac 
{{\it s_2}\,{\it s_3}\,{\it s_4}\,{\it s_6}}{144}}+\frac{{\it s_{15}}}{15}+{\frac {{{
\it s_5}}^{3}}{750}}-\frac{1}{26}\,{\it s_2}\,{\it s_{13}}+{\frac {{{\it s_3}}^{5}}{
29160}}+{\frac {{\it s_{10}}\,{\it s_2}\,{\it s_3}}{60}}+{\frac {{\it s_2}\,
{\it s_4}\,{\it s_9}}{72}},
\end{eqnarray*}
and gives a non zero answer
\[
{\mathcal B}=
{\frac {9639\,{c}^{15}}{17592186044416}}.
\]

\section{Nonlinear algebraic conditions}
\label{section4}
We now move to the second source of obstructions, namely the
nonlinear condition $J^2=-\mbox{Id}$. In the index notation this is
\be
\label{alg_1}
{\omega^a}_b{\omega^b}_c+\frac{1}{n}|\omega|^2{\delta^a}_c=0.
\ee
The case $n=4$ was treated in \cite{DT}, so in the Theorem below we shall assume that $n>4$.
\begin{theo}
\label{prop2}
There is a one--to--one correspondence between K{\"a}hler metrics in a
conformal class and parallel sections $\Psi$ of the vector bundle $(E,
\mathcal{D})$ from Theorem \ref{theo01} such that \be
\mathcal{Q}(\Psi)=0, \ee where $\mathcal{Q}$ is the set of non-linear
algebraic conditions given by (\ref{alg_1}) and
\begin{eqnarray}
&&\mu_{abc}+\frac{3n}{|\omega|^2}\omega_{[ab}{\omega^d}_{c]}K_d=0,
\label{mufromK}\\
&&\Sigma_{ab}=-\Big(\frac{n}{2}|\omega|^{-2}|K|^2+\frac{n}{4(n-2)(n-4)} C_{cdef}\omega^{cd}\omega^{ef}\Big)\omega_{ab}
\label{formula_for_sigma}\\
&&+\frac{1}{2(n-4)}C_{cdab}\omega^{cd}+
2n|\omega|^{-2}K_c{\omega^{c}}_{[b}K_{a]} \quad (n>4).\nonumber
\end{eqnarray}
\end{theo}
We shall split the proof into two steps. 
\begin{prop}
\label{prop21}
Solutions $\omega_{bc}$ to the CKY equation (\ref{mdprol0}) that also satisfy (\ref{alg_1})
correspond to conformally related K\"{a}hler metrics iff  (\ref{mufromK}) holds.
\end{prop}
\noindent
    {\bf Proof.} $\Rightarrow :$ In this direction the  result is immediate
    as  \eqref{solution_1} implies \eqref{mufromK}. 

    \medskip

  \noindent   $\Leftarrow :$   Assume (\ref{mdprol0}) and (\ref{alg_1}) hold. We would
  like to deduce (\ref{solution_1}), as then a conformal factor can be
  found which turns $g$ into a K\"ahler metric.
 Differentiating the condition (\ref{alg_1}) leads to
$$ \omega_{a}{}^b\nabla_d\omega_{bc} + \omega_{c}{}^b\nabla_d\omega_{ba}
=2\Omega^{-2}g_{ac}\Upsilon_d,
$$
where we defined the positive function $\Omega$ by $|w|^2_g=n\Omega^{-2}$, and as usual 
$\Upsilon_a : =\Omega^{-1}\nabla_a\Omega$.
Substituting (\ref{mdprol0}) yields
\begin{equation}\label{key}
 \omega_{a}{}^b\mu_{dbc} + \omega_{c}{}^b \mu_{dba} +
\omega_{ad}K_c +\omega_{cd}K_a
-g_{dc}\omega_a{}^b K_b-g_{da}\omega_c{}^bK_b
=2\Omega^{-2}g_{ac}\Upsilon_d.
\end{equation}
From this we find easily the second part of (\ref{solution_1}), as follows. Contracting 
(\ref{key}) with $g^{ac}$  and $g^{dc}$ respectively, and taking an appropriate linear combination
 of the resulting two equations yields
\be
\label{KKmm}
(n-2) K_c\omega^c{}_a + \omega^{bc}\mu_{abc}=0,
\ee
which then implies 
\begin{equation}\label{rho}
K_a=\omega_{ab}\Upsilon^b.
\end{equation}
Thus (\ref{key}) now gives a stronger relation
between $K$ and $\mu$, namely
$$
\omega_{a}{}^b\mu_{dbc} + \omega_{c}{}^b \mu_{dba} +
\omega_{ad}K_c +\omega_{cd}K_a
-g_{dc}\omega_a{}^bK_b-g_{da}\omega_c{}^bK_b
=-2g_{ac}\omega_d{}^bK_b ,
$$
or equivalently in terms of $\Upsilon$
\begin{equation}\label{key2}
\omega_{a}{}^b\mu_{dbc} + \omega_{c}{}^b \mu_{dba} +
3\omega_{a}{}^b\Upsilon_{[d}\omega_{bc]} + 3\omega_{c}{}^b\Upsilon_{[d}\omega_{ba]}=0 .
\end{equation}
Now clearly $\mu_{abc}=-3\Upsilon_{[a}\omega_{bc]}$ is a solution of
(\ref{key2}). Note however that, by linearlity, (\ref{key2}) only
determines $\mu_{abc}$ up to the addition of three--forms which belong
to the kernel of a linear operator $T:\Lambda^3\rightarrow
(\Lambda^1\odot\Lambda^1)\otimes\Lambda^1$ given by
\begin{equation}\label{key3}
T(\tau)_{abc}={{\omega_{a}}^d} \tau_{bcd}+ {{\omega_{b}}^d} \tau_{acd}.
\end{equation}
We can decompose
$\Lambda^3$ orthogonally  into parts that are trace-free and pure trace with respect to $\omega$:
$$
\Lambda^3 = \mathring{\Lambda}^3\oplus \Lambda^{1, 3} .
$$
The kernel of $T$ consists of 3-forms that are in $\Lambda^{3, 0}\oplus
\Lambda^{0,3}$ with respect to $J$, but all we need to know currently
is that if $\tau\in\mbox{Ker}\;(T)$ then $\omega^{ab}\tau_{abc}=0$ (as follows immediately by contracting (\ref{key3}) with $g^{ab}$), so
$\ker (T)\subset\mathring{\Lambda}^3 $. 

Let us write $I$ for the linear sub-bundle of $(\Lambda^1\odot\Lambda^1)\otimes\Lambda^1$ consisting of elements of the form
$$
-3\omega_{a}{}^b\Upsilon_{[d}\omega_{bc]} - 3\omega_{c}{}^b\Upsilon_{[d}\omega_{ba]}.
$$   
Then $T^{-1}(I)$ is a linear sub-bundle of $\Lambda^3$. On the other hand, since
$\mu_{abc}=-3\Upsilon_{[a}\omega_{bc]}$ is a solution of (\ref{key2}),
it is clear that $\Lambda^{1,3}\subset T^{-1}(I)$ and every
element in $T^{-1}(I)$ is, pointwise, the vector sum of an element
$\Lambda^{1,3}$ with an element of $\ker(T)$. Putting these things
together it is clear that we obtain unique solutions $\mu$ to
(\ref{key2}) if we restrict to $\mu$ which are pure trace. Or in other
words $\mu$s such that their trace-free part is zero:
We can write the condition for $\mu$ explicitly 
as\footnote{Note that this holds identically if $n=4$, as any element of $\Lambda^3$ 
in four dimensions is of the form $\tau=*\kappa$ for $\kappa\in\Lambda^1$, and
\[
\omega^{ab}{\tau}_{abc}=\frac{1}{6}\omega^{ab}\varepsilon_{abcd} \kappa^d=\frac{1}{3}\omega_{cb}\kappa^b
\]
as $\omega$ is self--dual. The RHS of the expression above never vanishes unless $\kappa=0$.}
\be
\label{alg_2}
|\omega|^2\mu_{abc}-\frac{3n}{n-2}\omega_{[bc}\mu_{a]pq}\omega^{pq}=0.
\ee
If this condition holds together with (\ref{alg_1}) then, using \eqref{rho}, there exists a unique solution to (\ref{mdprol0}) given by (\ref{solution_1}). 
Using (\ref{KKmm}) we see that \eqref{alg_2} is equivalent to (\ref{mufromK}).
\koniec

\begin{lemma}
For $\mu$ and $K$ as defined in \eqref{mdprol0}, the following
identities hold \be
\label{identity_KK}
\mu_{abc}K^c=0, \quad \mbox{and}\quad {\omega^a}_{[b}\Sigma_{c]a}=0.
\ee
\end{lemma}
\noindent
{\bf Proof.} In the conformal to  K\"ahler scale we have $\mu_{abc}=-3\Upsilon_{[a}\omega_{bc]}$ and
$\Upsilon_a=n|\omega|^{-2}K^d\omega_{da}$ so that 
\[
K^c\mu_{abc}\sim (\Upsilon_a\Upsilon_b-\Upsilon_b\Upsilon_a)=0.
\]

For the second property, which shows that $\Sigma$ is Hermitian, we use that from the expression for $\Upsilon$  in
terms of $\omega$ and $K$, as above, we have
\[
\nabla_a\Upsilon_b=\Omega^2(2\Upsilon_a\Upsilon_b\Omega
+{\Rho_a}^d \omega_{cd}{\omega^c}_b+\Sigma_{ac}{\omega^c}_b-K^c\mu_{abc}-g_{ab}|K|^2+K_aK_b),
\]
and so 
\[
\nabla_{[a}\Upsilon_{b]}=-\Sigma_{c[a}{\omega^c}_{b]}.
\]
But of course $\Upsilon$ is exact, and hence closed.
\koniec
\begin{prop}
\label{prop22}
If $n>4$ then $\Sigma$ is determined by $\omega$ and $K$, and is given by
(\ref{formula_for_sigma}).
\end{prop}
\noindent
{\bf Proof.}
Consider (\ref{mdprol1}), (\ref{mdprol2}) and (\ref{KKmm})
and compute
\[
\nabla_a(\mu_{bcd}\omega^{cd})=(n-2)\nabla_a({\omega_b}^e K_e).
\]
The RHS gives
\be
\label{RHSs}
(n-2)(g_{ab}|K|^2-K_aK_b-\frac{1}{n}|\omega|^2\Rho_{ab}+{\omega_b}^e\Sigma_{ae})
\ee
and the LHS  gives
\be
\label{LHSs}
-g_{ab}\Sigma\cdot\omega-2\Sigma_{bc}{\omega^c}_a+{\mu_b}^{cd}\mu_{acd}+\frac{2-n}{n}|\omega|^2\Rho_{ab}-\frac{1}{2}{C_{cda}}^p\omega_{pb}\omega^{cd}
\ee
where we have used (\ref{identity_KK}).

Next use
\[
\mu=3\Upsilon_{[a}\omega_{bc]}, \quad \mbox{and} \quad
\Upsilon_a=n|\omega|^{-2}K_b{\omega^b}_c
\]
to compute
\[
{\mu_b}^{cd}\mu_{acd}=2|K|^2g_{ab}-2K_aK_b+(n-4)n |\omega|^{-2} K_pK_q{\omega^p}_a
{\omega^q}_b.
\]
Substituting this formula in (\ref{LHSs}), comparing with (\ref{RHSs}) 
and contracting both sides with ${\omega^b}_e$
gives
\be
\label{almost_there}
\frac{n-4}{n}|\omega|^2\Sigma_{ae}=-(\Sigma\cdot\omega)\omega_{ae}+\frac{1}{2n}|\omega|^2C_{cdae}\omega^{cd}+(4-n)|K|^2\omega_{ae}+2(n-4)K_b
     {\omega^b}_{[e}K_{a]}.  \ee
Contracting this with $\omega^{ae}$ gives an expression for
$\Sigma\cdot\omega$ which we can substitute back to
(\ref{almost_there}).
This gives (\ref{formula_for_sigma}).  
\koniec
\noindent
Theorem \ref{prop2} now follows from 
Proposition \ref{prop21} and Proposition \ref{prop22}. 
\vskip5pt

The non-linear conditions  in Theorem \ref{prop2} trace a variety $\mathcal{S}$
in the fibres of the prolongation bundle $E$. If $n>4$ then 
\be
\label{constraint_surface}
\dim{\mathcal{S}}\leq \frac{1}{4}(n^2+2n+4).
\ee
To see it, note that  in Theorem \ref{prop2} both  $\Sigma$ and $\mu$ have been determined in
terms of $\omega$ and $K$. Substituting 
the expression for $\Sigma_{ab}$ into the expression (\ref{mdprol3})  for $\nabla_a\Sigma_{bc}$
could lead to an algebraic condition only involving $K$ and $\omega$. We claim that, at least in the conformally flat case, this condition is an identity, and does not lead to any further constraints on $K$. Indeed, we can choose a flat metric $g$ in the conformal class in which case 
(\ref{formula_for_sigma})  and  (\ref{mdprol3})
reduce to 
\[
\Sigma_{ab}=-\frac{n}{2}|\omega|^{-2}|K|^2\omega_{ab}
+ 2n|\omega|^{-2}K_c{\omega^{c}}_{[b}K_{a]}, \quad \nabla_{a}\Sigma_{bc}=0.
\]
Substituting the first expression into the second, and eliminating the
derivatives of $\mu, \omega$ and $K$ using the prolongation connection
leads to an identity. Therefore we can specify the $n$ components of
$K$ which are unconstrained, and the components of skew form $\omega$
which squares to a pure trace.  To count those set $n=2m$, take
\[
J=\left(\begin{array}{cc}
0 & I_m\\
-I_m& 0
\end{array}\right).
\]
and consider 
\[
\omega=J+\epsilon A, \quad\mbox{where}\quad \mbox A=\left(\begin{array}{cc}
a & b\\
c& d
\end{array}\right).
\]
Impose $\omega^2=-I_{2n}$. This gives $A=JAJ$, and
gives $b=c, a=-d$. If $\omega$ is a two--form
then $a=-a^T$, and only the skew--part of $b$ contributes so we can take $b=-b^T$
which gives a total of $m(m-1)$ components. There is one remaining component corresponding to the choice of an overall scale of $\omega$. Therefore the dimension of $\mathcal{S}$ is
$m(m-1)+1+2m$ which gives (\ref{constraint_surface}).
\vskip5pt
In dimension four, where the constraints on $(\Sigma, \omega)$ have been solved using self-duality \cite{DT},
and the bundle $E$ has been identified with the the bundle of self-dual
tractor three--forms which has rank $10$. The tractor approach for the general $n$ will be discussed in the next Section.
\section{Tractors}
\label{sectractors}
\def\frak{\mathfrak}
\def\Bbb{\mathbb}
\def\Cal{\mathcal}
\let\phi\varphi
\newcommand{\intprod}{\mathbin{\raisebox{\depth}{\scalebox{1}[-1]{$\lnot$}}}}

\newcommand{\cF}{\mathcal{F}}
\newcommand{\x}{\times}
\renewcommand{\o}{\circ}
\newcommand{\into}{\hookrightarrow}
\newcommand{\al}{\alpha}
\newcommand{\bet}{\beta}
\newcommand{\ga}{\gamma}
\newcommand{\de}{\delta}
\newcommand{\epsi}{\epsilon}
\newcommand{\ka}{\kappa}
\newcommand{\la}{\lambda}
\newcommand{\ph}{\phi}
\newcommand{\ps}{\psi}
\renewcommand{\th}{\theta}
\newcommand{\si}{\sigma}
\newcommand{\ze}{\zeta}
\newcommand{\Ga}{\Gamma}
\newcommand{\La}{\Lambda}
\newcommand{\Ph}{\Phi}
\newcommand{\Ps}{\Psi}

\newcommand{\Si}{\Sigma}
\newcommand{\Up}{\Upsilon}
\newcommand{\Ups}{\Upsilon}
\def\Rho{\mbox{\textsf{P}}}

\newcommand{\bg}{\boldsymbol{g}}

\newcommand{\barm}{\overline{M}}
\newcommand{\vol}{\operatorname{vol}}
\newcommand{\im}{\operatorname{im}}
\newcommand{\id}{\operatorname{id}}
\newcommand{\Fl}{\operatorname{Fl}}
\newcommand{\Ric}{\operatorname{Ric}}
\newcommand{\End}{\operatorname{End}}
\newcommand{\nrho}{\overset{\rho}{\nabla}}
\newcommand{\ntrho}{\overset{\tilde\rho}{\nabla}}
\newcommand{\nhrho}{\overset{\hat\rho}{\nabla}}

\newcommand{\rpl}                         
{\mbox{$
\begin{picture}(12.7,8)(-.5,-1)
\put(0,0.2){$+$}
\put(4.4,3.1){\oval(8,8)[r]}
\end{picture}$}}

\def\si{\sigma}
\def\W{\mathbb{W}}
\def\X{\mathbb{X}}
\def\Y{\mathbb{Y}}
\def\Z{\mathbb{Z}}
\newcommand{\lpl}{
  \mbox{$
  \begin{picture}(12.7,8)(-.5,-1)
  \put(2,0.2){$+$}
  \put(6.2,2.8){\oval(8,8)[l]}
  \end{picture}$}}
\def\form#1{\mathbf{#1}}
\def\dform#1{\dot{\mathbf{#1}}}
\def\ddform#1{\ddot{\mathbf{#1}}}
\def\dddform#1{\dddot{\mathbf{#1}}}
\def\ph{\varphi}
\def\rh{\rho}
\def\Up{\Upsilon}
\def\na{\nabla}
\def\de{\delta}
\newcommand{\ce}{{\Cal E}}
\newcommand{\nn}[1]{(\ref{#1})}
\def\Cal{\mathcal}
\newcommand{\bc}{\boldsymbol{c}}

The aim of this Section is to outline how the prolongation of the conformal Killing-Yano equations 
in Theorem \ref{theo01}
and the associated non--linear conditions on the parallel sections (Theorem \ref{prop2}) can be formulated  in the tractor language of \cite{BEG}.

In this section by a conformal manifold we mean a manifold equipped with an
equivalence class $\bc$ of Riemannian metrics such that if
$g,\hat{g}\in \bc$ then $ \hat{g}=\Omega^2 g$ for some positive
function $\Omega$. With $\cF$ denoting the frame bundle, the bundle of densities of
weight $w\in \mathbb{R}$ is the associated line bundle
$\ce[w]=\mathcal{\cF}\times_{\rho_w} \mathbb{R}$ where $\rho_w$ is the
representation of $GL(n,\mathbb{R})$ on $\mathbb{R}$ given by
$$
\mathbb{R}\times GL(n,\mathbb{R})\ni (x,A) \mapsto |\det (A)|^{-w/n} x\in \mathbb{R}.
$$ 
Therefore there exists a canonical isomorphism $\ce[2n]\cong \otimes^2(\Lambda^{n}TM)$. 
Using the conformal
structure a section $\phi\in \Gamma(\ce[w])$ can be identified with an
equivalence class of metric function pairs $[(g,f)]$ where
$$
(\Omega^2 g, \Omega^w f)\sim (g,f) .
$$
For any weight $w$, the bundle $\ce[w]$ is oriented and we write
$\ce_+[w]$ for the positive elements.

Using this notation it is straightforward to see that the conformal
structure determines a tautological section $\bg$ of $S^2T^*M\otimes
\ce[2]$ that we term the conformal metric. Then the metrics in a conformal
class correspond 1-1 with sections  $\si\in\Gamma( \ce_+[1])$ by the formula
$$
g=\si^{-2} \bg.
$$
We call $\si \in\Gamma( \ce_+[1]) $ (or the corresponding $g\in \bc$) a
{\em scale}. In the following we mainly follow \cite{BEG} and \cite{CG}.

\subsection{The standard tractor bundle and normal tractor connection}\label{basict}

On a conformal manifold $(M, \bc)$ there is, in general, no preferred connection on $TM$ but (in dimensions $n\geq 3$) there exists a connection  on a rank-2 extension
    \begin{equation}\label{tsplit}
    \cT=\ce[1]\lpl  T^*M[1]\lpl \ce[-1] ,
    \end{equation}
    that we call the {\em conformal tractor bundle} $\cT$. The right
    hand side of \eqref{tsplit} gives the composition series of $\cT$
    and means that $\ce[-1]$ is a canonical sub-bundle of $\cT$. The
    quotient of $\cT$ by this has $T^*M[1]$ as a canonical sub-bundle,
    and then the quotient by this is $\ce[1]$.  The tractor bundle
    $\cT$ will be denoted $\ce_A$ in the abstract index notation.

 Given $g\in
 \bc$ the tractor bundle splits into a direct sum
$$
\cT\stackrel{g}{=}\ce[1]\oplus T^*M[1]\oplus \ce[-1].
$$
So $V_B\in \Gamma(\ce_B)$ can be then represented  by a triple
$$
V_B\stackrel{g}=\begin{pmatrix}\si \\
\mu_b \\\rho \end{pmatrix} .
$$
The {\em normal tractor connection} is
$$
\nabla^\cT_a(\si,\mu_b,\rho )= (\nabla_a\si-\mu_a,~\nabla\mu_b+P_{ab}\si +\bg_{ab}\rho,~\nabla_a\rho-P_{ab}\mu^b ).
$$
This connection acts on tensor powers of $\cT$, and $\nabla^\cT$
preserves a conformally invariant {\em tractor metric} $h$ that is given as a
quadratic form by
$$
\cT\ni V=(\si,\mu_b,\rho )\mapsto 2\si \rho +\mu_b\mu^b=h(V,V). 
$$ In abstract indices we denote this $h_{AB}$ and use it to raise and
lower tractor indices.
It is convenient to introduce the  algebraic splitting operators $X^A,Y^A, Z^{Ab}$ that encode the slots,
\begin{equation}\label{XYZ}
  V^A= \si Y^A + Z^{Ab}\mu_b+X^A\rho,
\end{equation}
and we will need this below.
For example in the slot notation $X_B$ is represented by $(0,~0,~1)$.

There is a conformally invariant differential
operator $D:\Gamma(\ce[1])\to \cT$ given, in a scale $g$, by
$$
\si\mapsto D_B\si\stackrel{g}{=}(\si,~\nabla_b\si,-\frac{1}{n}(\Delta\si+ P\si) ),
$$
where $P:=\bg^{ab}P_{ab}$. Given the splittings as described this
is determined by the tractor connection formula. (Alternatively, when
the tractor bundle is constructed via jets, this operator actually
determines the splitting of tractor bundle into the triples
\cite{CG}.)  It is termed a splitting operator as the composition $X^B
D_B$ is the identity on $\Gamma(\ce[1])$.  If $\si$ is a scale, then
we define
$$
I_B:= D_B\si ,
$$ and call $I_B$ the corresponding {\em scale tractor}. So $\si=X^A
I_A$.
The squared length of the scale tractor recovers a multiple of
the scalar curvature of the metric $g=\si^{-2}\bg$:
$$
I^AI_A=h^{AB}I_AI_B=-\frac{1}{n(n-1)}R.
$$

One reason that $D_B$ is important is that if $V_B\in \Gamma(\ce_B)$
is parallel for the tractor connection then $V_B=D_B\tau$ for some
$\tau \in \Gamma(\ce[1])$.

\subsection{3-form tractors}
It follows 
from the semi-direct composition series  of $\cT$ that  the
corresponding decomposition of $\Lambda^3 \cT$ is
\begin{equation} \label{comp_series_form}
  \mathcal{E}_{[ABC]} = 
  \mathcal{E}^{2}[3] \lpl \left( \mathcal{E}^3[3] \oplus
  \mathcal{E}^{1}[1] \right) \lpl \mathcal{E}^{2}[1],
\end{equation}
where $\ce^k[w]$ denotes  $\Lambda^k(T^*M)\otimes \ce[w]$.

3-form tractors are useful for studying the conformal Killing Yano equation.
For $\si_{ab}\in\Gamma(\ce_{[ab]}[3])$ let us write
$$
KY(\si)_{abc}: =
 \nabla_{a}\si_{bc}  -\nabla_{[a}\si_{bc]} +\frac{2}{n-1}\bg_{a[b}\nabla^p\si_{c]p} \, .
 $$
 This is conformally invariant and $\si$ is a conformal Killing-Yano tensor if
 \begin{equation}\label{cKYeq}
KY(\si)_{abc} =0.
 \end{equation}

Now, a choice of metric $g$ from the conformal class
determines a splitting of the bundle $\Lambda^3 \cT$ into four components (a
replacement of the $\lpl$s with $\oplus$s is effected) so that a 3-tractor
$\Phi$ can be written a 4-tuple
$$
\Phi_{ABC}\stackrel{g}{=} \begin{pmatrix}
    \si_{bc} \\ \nu_{abc}\quad \ph_{c} 
      \\ \rh_{bc}
  \end{pmatrix}
$$
where $\si_{bc}\in \ce^2[3]$, $\nu\in\mathcal{E}^{3}[3] $ and so forth.

Given $\si_{bc}\in \Gamma(\ce^2[3])$ there is a conformally invariant
differential splitting operator
$$
L: \Gamma(\ce^2[3]) \to \Gamma (\Lambda^3 \mathcal{T}) ,
$$
determined by the tractor connection, and given by
\begin{equation}\label{promap}
\Gamma(\ce^2[3])\ni \si_{bc}\mapsto L(\si)\stackrel{g}{=} \begin{pmatrix}
    \si_{bc} \\ \na_{[a}\si_{bc]} \qquad  \frac{2}{n-1} \na^b\si_{bc} 
      \\  \frac{1}{2n} \na^p KY(\si)_{pbc} 
      - \frac{1}{n-1} \na_{b} \na^p \si_{p c} 
      - P_{b}^{\ p} \si_{p c}^{} \Bigr) 
  \end{pmatrix}\in\Gamma (\Lambda^3 \mathcal{T}) ,
\end{equation}
see \cite{GS}.
Now the key importance of $L$ is that is is related to
the prolongation connection $\mathcal{D}$ of Theorem \ref{theo01} for
the conformal Killing-Yano equation \eqref{cKYeq}. The following is a
special case of Theorem 3.9 in \cite{GS}.
\begin{prop}
  There is a conformally invariant connection  
  $\mathcal{D}$  on $\Lambda^3\cT$
  with the property that
  \begin{enumerate}
  \item
  $
{\mathcal D}\Phi=0 
  $
implies that
$$
\Phi=L(\si_{ab}) \quad\mbox{and} \quad KY(\si_{ab})=0 ;
$$
\item If $KY(\si_{ab})=0$ then 
$\mathcal{D}(L(\si_{ab}))=0 $. 
\end{enumerate}
This has the form
  $$
{\mathcal{D}}_a \Phi_{BCD}=\nabla_a \Phi_{BCD} + (\kappa\sharp \Phi)_{aBCD}
$$
where $\nabla_a$ is the normal tractor connection and $\kappa\sharp \Phi$ a linear action of its curvature on $\Phi$.
\end{prop}
\noindent The details of $(\kappa\sharp \Phi)_{aBCD}$ will not be needed below, but, with a little translation, they can
be read-off from (\ref{mdprol3}).

We want to apply this to $\om_{ab}\in \Gamma(\ce^2[3])$. If we assume that
$\om_{ab}$ satisfies the conformal Killing-Yano equation then the image of
\eqref{promap} simplifies to
\begin{equation}\label{ssplit}
L(\om_{ab}):=  \begin{pmatrix}
    \om_{bc} \\ \na_{[a}\om_{bc]} \qquad  \frac{2}{n-1} \na^b\om_{bc} 
      \\        - \frac{1}{n-1} \na_{b} \na^p \om_{p c} 
      - P_{b}^{\ p} \om_{p c}^{}
  \end{pmatrix}.
\end{equation}

\subsection{Characterisations of conformally K\"{a}hler}

For $\om_{ab}\in \Gamma(\ce^2[3])$ let us write
\begin{equation}\label{sieq}
  \si:=|\om|:=\sqrt{\frac{1}{n}\om_{ab}\om^{ab}},
\end{equation}
where indices have been raised
by the conformal metric, so if $\om$ is non-zero then
$$
\si\in \Gamma(\ce_+[1])
$$
is a distinguished scale determined by $\om$.

\begin{prop} \label{trchar1}
  The conformal class $\bc$ contains a K\"ahler metric iff  there exists $\omega\in
  \Lambda^2(M)$ such that
\begin{equation}\label{herm1}
\om^a{}_{c}\om^{c}{}_b=-\si^2\delta^a_{b},
\end{equation}
and
\begin{equation}\label{char22}
X\wedge I \wedge \Phi=0,
\end{equation}
 where  $\Phi=L( \om)$ and
 $I:= D\si$, with $\si$ defined by \eqref{sieq}.
\end{prop}
\noindent
{\bf Proof.}  $\Leftarrow$: From \eqref{herm1} we have that $\si$ is a
scale and that $\si^{-3}\om_{ab}$ is Hermitian for the metric
$g:=\si^{-2}\bg$. The Levi-Civita for $g$ preserves $\si$,
i.e. $\nabla^g\si=0$. Working in this scale we have 
  $$
  D_A\si=\si Y_A -\frac{\si}{n} P X_A .
  $$
  Thus
  $$
X\wedge I \wedge \Phi= \si X\wedge Y  \wedge \Phi .
$$
So, with $\Phi=L(\om)$, \eqref{char22} exactly captures the condition that 
the $Z\wedge Z\wedge Z$-slot of \eqref{ssplit} is zero, that is that $\si^{-3}\om$ is closed. Thus $\si^{-3}\om$ is the K\"ahler form for the K\"ahler metric $g$.

$\Rightarrow$: If $\si^{-3}\om$ is a K\"ahler form for a metric $g=\si^{-2}\bg$ then we have \eqref{herm1}. Moreover
the Levi-Civita connection of $g$ preserves $\si$ and $\si^{-3}\om$. Thus, in particular,  the latter is closed and co-closed in the scale $g$, and
$L(\om)$ takes the form
\begin{equation}\label{zz}
L(\om)\stackrel{g}{=}\begin{pmatrix}
    \om_{bc} \\ 0 \qquad  0 
      \\      
      - P_{b}^{\ p} \om_{p c}^{}  .
  \end{pmatrix}
\end{equation}
From this  it is evident that \eqref{char22} holds.
  \koniec
\noindent
{\bf Remark.}
 From the last display we see that for $\Phi=L(\om)$ satisfying
  \eqref{herm1} and \eqref{char22} we must also have
  \begin{equation}
X\hook I\hook \Phi =0.
  \end{equation}
  This is the co-closed condition.

Next,  from the same display we also see that in the case that \eqref{herm1} and \eqref{char22} hold, then the squared length of $L(\om)$, i.e.
  $$
\Phi^{ABC}\Phi_{ABC}
$$
is a non-zero constant times the scalar curvature of the K\"ahler metric $g=\si^{-2}\bg$.  
\vskip5pt

It is well known that on a conformal structure a metric $g$ is
Einstein iff there is a parallel tractor $I_A$ and $g=\si^{-2}\bg$
where $\si=X^A I_A\in \Gamma (\ce_+[1])$ is nowhere zero. If a tractor $I_A$ is parallel for the normal tractor connection then $I_A=D_A\si$ for some $\si\in \Gamma(\ce[1])$. These results follows
from the construction of the tractor connection in \cite{BEG}, as 
discussed in \cite{G,GNur-obstructions}.

Thus one immediately has the the following result.
\begin{prop} \label{trchar2}
  The conformal class $\bc$ contains a K\"ahler--Einstein metric if
  there exists $\omega\in \Lambda^2(M)$ such that
$$
\om^a{}_{c}\om^{c}{}_b=-\si^2\delta^a_{b},
$$
and
$$
X\wedge I \wedge \Phi=0,
$$
 where
 $I=D\si$ is parallel  for the normal tractor connection and  $\Phi=L( \om)$.
\end{prop}

A special case is when the Einstein structure considered is Ricci
flat. The length of the scale tractor is a multiple of the scalar
curvature. Thus $g=\si^{-2}\bg$ is scalar flat iff $I:=D\si$ is
null. On the other hand from \eqref{zz} we see that if the K\"ahler scale is Ricci flat then $I\wedge L(\om)=0$. So we have the following result.

\begin{prop} \label{trchar3}
The conformal class $\bc$ contains a Ricci--flat K\"ahler metric iff there exists $\omega\in \Lambda^2(M)$ such that
$$
\om^a{}_{c}\om^{c}{}_b=-\si^2\delta^a_{b},
 $$
and
$$
I \wedge \Phi=0,
$$
 where
 $I=D\si$  is parallel and null for the normal tractor connection, and $\Phi=L(\omega)$.
\end{prop}

\begin{remark} Since ${\mathcal D}\Phi=0$ implies $\Phi=L(\om)$ for some $\omega\in \Gamma(\ce^2[3])$, there are (slightly weaker) variants of these propositions where we replace, in each case, the condition $\Phi=L(\om)$ with  $\mathcal{D}\Phi=0$.
  \end{remark}

\section{Outlook} We have constructed a rank $n(n+1)(n+2)/6$ vector bundle $E$ with a connection ${\mathcal D}$  over a Riemannian manifold $(M, g)$ of even dimension $n$, such that 
the ${\mathcal D}$--parallel sections of $E$ belonging to a certain
non--linear variety ${\mathcal S}$ in the fibres of $E$ are in
one-to-one correspondence with K\"ahler metrics in a conformal class
of $[g]$. The construction of the connection followed from the
prolongation of the conformal Killing--Yano (CKY) tensor equation
\cite{Sem, GS, DT}, and the construction of ${\mathcal S}$ resulted
from exploring the differential consequences of $J^2=-\mbox{Id}$,
where the endomorphism $J:TM\rightarrow TM$ is the complex structure
of the K\"ahler form.

The integrability conditions for the existence of the parallel
sections in ${\mathcal S}$ imply that the conformal Weyl tensor of $g$
is of the algebraic type--$D$. We have established an explicit
algebraic obstruction for this which makes the results relevant in
general relativity of type--$D$ spaces in dimension higher than four
\cite{coley, pravda, mason}.

The conformal Killing--Yano tensors which underlie our work give rise
to hidden symmetries of gravitational instantons \cite{jezierski, DT,
  DT2, houri, araneda}, as well as to first integrals of the conformal
geodesics \cite{AR, DT3}.  The obstructions we have constructed can be
of separate interest in deciding whether a given metric (Lorentzian or
Riemannian) admits such hidden symmetries, or whether a conformal
geodesic motion is integrable.

Finally, there is a connection with the tractor approach to conformal
differential geometry \cite{BEG}: the prolongation bundle $E$ in our
work can be identified with a parallel transport condition on
$\Lambda^3({\cT})$, where ${\cT}\rightarrow M$ is the rank--$(n+2)$
tractor bundle. It is however the case that the connection induced on
$\Lambda^3({\cT})$ by the standard tractor connection on ${\cT}$
differs from the prolongation connection ${\mathcal D}$ we have
constructed on $E$ in Theorem \ref{theo01}. It would be interesting to
reformulate the non--linear algebraic conditions on the parallel
sections of $E$ in our Theorem \ref{prop2} purely in terms of
tractors. This is essentially implicit in Proposition \ref{trchar1} as
$D\si$ and $\om$ can each be expressed algebraically in terms of
$\Phi$ (as the prolonged system is closed). However it would be useful
to find a simpler and explicit description.  In \cite{DT} this problem
has been solved in dimension $n=4$, where the non--linear conditions
reduce the bundle $E$ to the rank--10 bundle ${\Lambda^3}_+({\cT})$ of
self--dual tractor three--forms. The problem of finding an analogue of
this remains open for $n>4$.

\end{document}